\documentclass[11pt]{amsart}
\usepackage{a4,palatino,mathplm,txmac}
\usepackage[breaklinks]{hyperref}
\makeatletter

\advance\oddsidemargin by -1.25cm
\advance\evensidemargin by -1.25cm
\advance\textwidth by 2.5cm
\advance\footskip by 3.5mm
\advance\headsep by 2.5mm
\advance\voffset by 11mm

\begin{document}
\title{Conduct and Correctness in Mathematical Publishing}

\author{Alexander Stoimenow}
\let\@br\\
\address{Department of Mathematics, Keimyung University, \@br
\indent Darseo-Gu, Dalgubeoldaero 2800, Daegu 704-701, Korea}%
\email{stoimeno@stoimenov.net}%
\urladdr{http://stoimenov.net/stoimeno/homepage/}%
\let\@cite@ofmt\relax 
%allows reference labels to be treated as usual text and line-broken

\begin{abstract}
This is an essay in which I try to express my fear about the
establishment of a culture of publishing where no one is willing
to take responsibility for the correctness of mathematics, and
readers finding mistakes in published proofs are stamped as
outcasts, because they are deemed to target the reputation of
authors and journals. \\[2mm]
{\em keywords:}\ Publication of research, correctness, peer
  review, conflict of interest
\end{abstract}

\def\dd{\discretionary{}{}{}}

\maketitle

\section{Introduction: journal publishing for mathematical research}

Research belongs, to a larger or lesser extent, to the life
of every academic mathematician. Unlike other sciences, the
validity of mathematical research depends on calculations
and logical steps which derive a new result (theorem) from
previously known ones. Nowadays mathematics has grown more
complex than the expertise of individual mathematicians, and
the knowledge needed to understand and apply a theorem often
differs considerably from the one needed to understand its
proof. However, if a flaw is found in the proof of a theorem
(unless other correct proofs exist), everything building on
that `theorem' is potentially also flawed. On the other hand, 
electronic media have facilitated the communication of research,
and therewith also of poor one. This dramatically increases the
importance of journals as places to look for reliable material. 

The process of publishing in journals briefly goes thus. An author
writes up the results of his research in a paper, and chooses a
journal to submit. The editor of that journal examines, with the
help of one or more referees, the suitability (including correctness)
of the content, and decides whether to publish it or not. In latter
case, the author may submit the paper to another journal, where the
procedure repeats. (See \cite{Bernstein,Bernstein2} for a general
code and \cite{B} for an interdisciplinary survey of peer review.)

As the process has been arranged in mathematics, the referee often
knows who the author is, but not vice-versa. This practice is
justified by assuming that a referee (who is thought to be a mature,
mathematically and ethically responsible individual) would not
personalize this relationship. Still there is, in theory, no
restriction on what kind of comment a referee can make on a paper,
nor on which of the material an editor receives he chooses to
publish and which not. In
the end, authors submitting papers to a journal and its readers have
not more than a word of good faith of the editors that the material
is properly evaluated and selected. While this may generally be
the case, and there is always some margin of personal preference
(and human error) in the process, it should be unnecessary to
point out that the opportunity to abuse such a system is virtually
unlimited. As an author, I would hope to, and often do, respect the
work of editors (among which I am not) and referees (of which I am
one occasionally). Still, all too often I am left wondering which
scientific principle can explain what I see. 

Ultimately, when mathematical correctness issues got involved, my
worries became serious enough to motivate the following report. It
centers around a concrete case with a reputed American research
monograph publication M.
Before turning to this case, however, I chronologically follow
several situations in its background. They link it with the issue
of journal publication in several ways and explain my relation to
the two authors of the monograph, whom I will call Prof.\ A and
Prof.\ B below\footnote{While I keep record of most (written)
correspondence from which I quote, I can, for suggestive reasons,
not afford to target directly diverse individuals and institutions.
I offer, though, to provide this documentation, to serve as a factual
basis independent of my interpretations.}.

So far little publicity has been offered in mathematics to this
fundamental issue~-- how to maintain correctness of the literature~--
and to its relation to journal organization. It is difficult to
say to what extent this status quo is owed to
lack of occasions where the issues surface, or to the amount of
neglect and suppression they experience. There is, however, some
evidence that my case is not isolated, while few opportunities are
established within mathematical community to debate such problems.
This is a reason for me to seek some venue for discussion here.

I am grateful to Rob Dickey for his valuable suggestions on a
previous version of this text.

\section{The interaction between referees and authors\label{SI}}

Mathematics has a smaller community than the one of most other
sciences. This situation can result in exceptional authority granted
to a few, as editors have a great deal of discretion and certain
individuals are frequently called upon to serve as a referee, perhaps
in part because of their high standing in the community, perhaps also
in part on their willingness to perform the role.

My experience with Prof.\ B exemplifies a result of this circumstance.
Our initially good relationship deteriorated when he refereed and
rejected a submitted paper of mine with the argument that a different
proof of my result was written in a (at that time) 7-year old
unpublished draft of his. While there might have been other meaningful
reasons to argue against publication of my submission, his own
failure to write up his results should not be considered as such.
I had received a copy of his draft (after submitting my paper)
and, without figures, found it unreadable. I then made a fatal
mistake of criticizing him directly (thereby revealing to him,
of course, that I knew he was the referee). I naively assumed
that a senior scientist, whom I also know personally, would be able
to react maturely in such a situation. I received no reasonable
response from him, but instead another pile of similarly looking
comments about my work from various other journals. After this
continued for a while, and I saw no other way to defend myself,
I put some of these reports on my website, seeking help to make
him stop his activity. Of course, I had again to attribute the
reports to him, without that an editor would confirm his identity
to me~-- it is a refereeing system that always protects the referee
at the cost of the author. (At the time our dispute broke out, B
had allegedly declared that he would never again referee my papers.)
Although I did (privately) hear some conciliatory voices, my attempts
to openly violate anonymity and object to various (also other)
referees' actions offended B, and not only him. An editor I pointed
to my website reacted as follows.
\begin{quote}
I am satisfied that the referees consulted were honest and did
professional jobs.

If your comments are criticisms of the refereeing process in general,
then you should find some other forum in which to debate them.

If they are specific criticisms of [our institution], and you genuinely
believe them, then you may do better to find some other journal to
submit your articles to. [$\ldots$]

I now consider this correspondence terminated.
\end{quote}

Is referee anonymity holy enough to teach young and unestablished
research\-ers to fear for their academic survival? B's "fight" has
been going on for about 13 years now, and includes a recent case where
I solved a problem formulated by himself, and had a coauthor who~--
like me~-- was continuously struggling for jobs, but~-- unlike me~--
was completely impartial to all this conflict. B's findings that
our paper had some weaknesses and the journal was too prestigious
might have been agreeable, but there was actually not one single
positive word in the report. My coauthor and I received then from
another journal a similar report, largely copied word by word from
the previous one. I know also of other referees accustomed to
such "copy-paste" practices, and I try hard to tolerate them as a
"feature" of the reviewing system. Formally, as long as only B's
comments on my papers are concerned, this can still go as his opinion
and, ultimately, there will be some editor who will not choose him
as a referee. Among others, the paper I mentioned he rejected at
the time appeared 5 years later, with additions and improvements,
in a much better journal. (In comparison, a completed, readable
version of his draft was published only recently, almost 20(!)
years after its announcement.) In some respects, we may say in this
case that the interaction with the peer review system `helped' me.
However, when the career of my collaborators is (also) derailed,
and definitely when readers of B's own papers are kept misled (as
I will describe in the monograph case in \S\ref{S4}), battling
me across editorial offices ceases to be a mere private matter.

Prof.\ A addressed the conflict like this: "He is now your enemy"
and "you cannot afford more enemies". First I tried to see this as
advice (also from others) about personal relationships. But with
A's involvement in the monograph case, I started questioning such a
position. Exercised indiscriminately, private hostilities can easily
turn into a device for transgressing whatever standards of conduct
exist~-- be it through intimidation, attacking at will, or ignoring
scientific disagreement. Who qualifies to be declared an "enemy"?
Those whose objections about one's work, or whose own work, one has
authority to spread doubt on? 

\section{Questions of editorial conduct\label{S3}}

Questionable practices also occur on the part of the editor.
Once, for example, I was sent to review a true (already published)
"masterpiece" for the online review database \cite{MSN} of the
American Mathematical Society (AMS). The paper had something
suggested as a "main result", whose proof was 8 lines, with this
being the only proof in the paper. The idea was, roughly speaking,
to observe that when $a\le b$, $b\le a$ and $b=c$, then $a=c$. I
did not find anything else really new inside. I perceive such an
article as a mockery, not only against myself~-- I had 9 papers
previously rejected by that journal, 8 without seeing any referee
report~-- but also against the efforts of dozens of other authors,
who submit much more serious material there, and many of whom are
"processed" in a similar way.

In the end, I decided I could
not write a review~-- I saw no way to write anything serious without
sounding offensive. Still, I could not refrain from sending an angry
letter to the editors, asking how a traditional top-30 journal could
accept something like this. (The managing editor was heartfully
thanked in the acknowledgment for his help to `improve' the paper.)
For various (personal) reasons, I felt unwilling to raise the issue
formally; the
editors, however, had little appreciation for that. As I had declared
that I might do otherwise in a following case, they referred to
these `threats' as a pretext to collectively refuse to handle my next
submission. So while the editors could appear to be comfortable with
their doings, it was an outsider who was creating disorder.

Prof.\ A was following up close my outrage about this paper, and
the treatment of mine. From his own experience, he describes a
possible editor's position as follows: "If his papers are good,
some journal would publish his papers", and "once I decided, it
is usually final". So, "[w]hat kind of response do you expect from
the editors? Apology?" Still, Prof.\ A\ admits that if a decision
is made "by a political reason", then "it is hard to prove it".
Thus, rather than plausibly addressing editorial difficulties or
following scientific guidelines, the opportunity is created to trump
up arbitrary ruling. (What proof should expose one? Criminalistic?
See also the quote of P.~Doty and the end of \S\ref{S8}.) However,
in response an editor might appeal to a universal (because
unobjectable) justification: difference in opinion. Indeed, any
perception of the quality or importance of a paper, no matter how
many would share it, is subjective~-- and hence can always be disputed
by a particular editor or referee. In the long run, if published
material is technically correct, there is no irreversible harm done
to mathematics. This is less certain, though, at the presence of
errors. Such a situation arose a few years later~-- and became the
occasion for my essay. 

\section{The monograph case\label{S4}}

During my research I discovered an error in the monograph of the
two authors. It was possible, to large extent, to repair the problem,
but it required some explanation. An 8-page note providing the
details was prepared and sent to the managing editor (in 6/07) to
present the issue to the authors. (I decided not to communicate
with the authors directly, after the explained experience, and
after having no reply from Prof.\ A\ to my previous email regarding
an error in another one of his papers.) I hoped for some serious
response~-- that either they could refute the objection, or write
an erratum based on the explanation. I received, though, just a very
brief comment from Prof.\ B, sent on behalf of both authors, a few
months later. He merely stated that "[w]e never claimed" what I
was referring to, and "[i]t seems to us that Stoimenow may not
understand details of our paper". It is true that the claim does
not appear explicitly, but it enters into the proof of one of
the main results. He also sent a 'revision' of this proof, but it
had corrected only a few typos. He did not address any of my concerns
(which suggested to me that he had not read my note at all, or wished
to quickly dismiss it), and I replied listing again one-by-one the
questions and objections. To this date I have not witnessed any further
mathematical comments from the authors addressing this proof.

First I had submitted, in 8/08, to the same publication M an own long
manuscript, in which I included the discussion of the error. Seven
months later, in 3/09, a message stated that M could not consider
the manuscript, partly due to backlog (of which I had been admittedly
advised in advance by the managing editor). Then I submitted my short
note as an erratum, assuming that the authors had not written their
own. (I was not told otherwise by the managing editor, and it
was meanwhile a year and 9 months since my initial indication of the
error.) After apparently discussing briefly with the authors, the
managing editor, however, declined to handle my correction, too (in
particular to have it examined by an independent party). I quote
his message:
\begin{quote}
Corrections are by the same authors. So this must be considered a
new article. It certainly is not appropriate for M by length
considerations. It certainly seems that you and the authors of the
original [monograph] do not agree on many things.
  
I think that if you want to publish this, you should submit to some
independent journal and have it refereed.
   
So we cannot consider this for M.
\end{quote}
He did not explain what issues I disagree on with the authors, and
did not ask my position on these issues. I received no response to a
further query, asking whether the authors submitted an erratum, and
in what form errata to M can be published. Later the same month I met
Prof.\ A at a conference, and he said that Prof.\ B would contact me
regarding the matter. Also, in 2/08, a third person had told me that
the authors plan to put a corrected version of the monograph on the
arXiv \cite{arXiv,Jackson}. So I waited. After a further ten months,
still nothing had happened.

Before I tried submitting the case elsewhere, a senior colleague (whose
similar experience I had found only weeks before \cite{Hill,Hill2})
advised me to contact the managing editor one more time directly. I 
repeated my questions to the editor and included the (major, and critical) 
parts of a draft of the present text that concerned his case. Despite being 
informed in the draft, among other things, of my conflict with Prof.\ B, 
he reaffirmed his determination not to handle (and ask independent opinion
on) my correction. He just forwarded to me again the old (and as I
explained, completely inadequate to me) response of Prof.\ B, and
declared that "the decision was not to submit an erratum to M".

Only at that point I started sending my correction, \emph{a fortiori},
to other journals. At the time of writing (06/12), after M, at least
8 other journals declined to consider, or rejected my erratum. Only
occasionally there are signs of a more careful
examination. From the editor of the fourth journal I received the
first, so-to-say, authorized word "that there is indeed a mistake
in the [A-B] paper". However, according to "a couple of experts",
"repairing that mistake is not such a big issue". I do not know
whether this, or what else, is a publicly agreeable judgement to
close the matter, after all disarray that occurred; for my reasons,
see the beginning of next section. The dilemma surfaces more
directly in one of the following referee reports:

\begin{quote}
This paper describes an error in a [monograph] that appeared
in [year]. Given that the [monograph] has been cited [X] times
(according to the MathSciNet count) it is not a major publication.
Nonetheless, it is of interest to the [specialized] community and
has generated further research. Thus a paper correcting an error
in it should be published. But does [your journal] publish papers
of that sort, or should there be an audience outside of [this
particular] community?
\end{quote}

In summary, 5 years after I originally raised my objections to M,
I still have seen no effort of others to address (or to allow me
to address) the issue officially.

\section{Mathematical lit(t)erature?\label{Sli}}

Flaws have been found also in diverse other papers, and that
some authors and editors are unenthusiastic in considering
these is hardly novel, also to myself. The reasons I attached
some significance to the problem with the monograph are that
\begin{list}{$*$}{\itemsep 3pt\relax\topsep 4pt\relax}
\item this is one of the main results of the monograph,
\item there is no (to-be-corrected) reprinting of research monographs
  scheduled (unlike e.g. undergraduate textbooks, to which one may
  tend to compare M rather than to research journals),
\item the error is subtle (it had not been noticed for about 15 years
  since the manuscript was written),
\item this is the only proof available,
\item other work (published also in other papers) depends essentially
  on this result, and
\item I know at least one other (published) paper which uses the
  argument and may (without that I went into details there) encounter
  the same problem.
\end{list}
Of course, I can unlikely win a
dispute about subjective views on importance of corrections. But while
it is suggestive that one shouldn't necessarily write an erratum for
every typo, it is hardly convincing when I am left to wonder whether 
my explanation of the error was carefully read and understood. There
is some consensus in the scholarly community that it is customary the
original authors to write an erratum, but in part we may also consider
this to be their (and not my) duty. 

What justifies such authors' conduct? Unlike possibly in tabloid
journalism, an argument in science that something (here, a correction)
is redundant to publish because presumably not of great interest
makes little sense: should we write no proofs because the fewest
read them? Do the authors reckon that I better not "[dis]agree on
many things" with them, when I "cannot afford more enemies"? Or do
they find in their great expertise confidence not to take seriously
my concern? Lang writes \cite[end of \S V.4]{Lang}: "I do not
recognize being an 'exceptional scientist' as a license to throw
one's weight around to avoid answering scientific criticisms."
I see his reasons. In the end, if littering around in one's
publications is a negotiable freedom, who will regulate what and
how to be traded for it, and who will organize the cleanup? Many
researchers know that finding a correct proof of a non-trivial
statement is often painstaking. How will one be motivated in going
through this when seeing others well received with such working
attitude?

Another point to keep in mind is that, even if a journal does not
remunerate its contributors, publications are not a charity enterprise.
There is thus little reason in expecting their correctness to be one,
instead of holding authors personally accountable for what they seek
employment, tenure or grants upon. However, as noted in \cite[V.4(b)]
{Lang}, in reality notorious cases throughout science (see, e.g., 
\cite{Bu,O+}) arise with established figures, to whom career and
funding pressures are, at least, not existential. Such individuals
know well that few of the affected can afford an open conflict and,
if a conflict indeed arises, can take advantage in it. Most people
naturally stay away from disputes, rather than investigating what
whose merits or stakes. Others go along with a scientist's prominence,
or at least are reluctant to openly and actively take position
against him (and whatever negative publicity is far more critical to
the less established ones). Someone wrote to me regarding my case:
"I just want to be friend with all involved parties." This is just,
unfortunately, not what readers seek in one's papers.

One other defense of one's conduct is to refer to the existence of
courts of law. Lang writes in \cite[\S V.3]{Lang} that "[s]uch a
point of view undermines the exercise of scientific responsibilities,
as distinguished from legal responsibilities." This way practices
like sloppy writing, refusal to correct errors, or lousing
around on junior colleagues, become acceptable, because they are
\em{not illegal}. Similar is an argument of the sort `this is
how it goes'. Accordingly, corruption in mathematics is legitimized
when a majority of mathematicians becomes (or because it already
is?) corrupt. (See \cite[p.447-451]{Lang2} for a high-profile
instance of such reasoning.) In many conversations Prof.\ A (and
not only he) has discouraged me from (and criticized others for)
settling scientific disputes legally. I find, though, increasingly
less reason to share his optimism about academic authorities'
self-policing. In result, I feel continuously pushed toward a
structure of behavior, where making a `career' means to let go
whatever higher ups please, in order to keep a chance of being
allowed one day to do similar with others.

One might try to imagine objective reasons for such an indefinite
delay, in that an author, as a human being, is subject to change of
duties in life. Much less so, however, is a journal (or monograph
series). This is why it is a plausible thought that it should be
ultimately the journal that has to properly address issues of the
correctness of the mathematics that has appeared in it, and to
provide ways for corrections to be published.

\section{The role of journals\label{JR}}

There seems no definite codex of editorial conduct, and so editors
have some freedom to set their own policies. But whatever these may
be, they should not oppose basic scientific principles.

Is it reasonable that potential rules regarding such matters as
authorship of errata, or adopting the authors' position in a dispute
with readers, are to take precedence over issues of correctness of
mathematics? If the journal cannot guarantee the author's responsible
involvement, it should grant the right to others to publish corrections,
and certainly have them carefully~-- and independently~-- considered.
For example, in the Hales-Hsiang case \cite{Hales,Szpiro}, the
\emph{International Journal} did publish Bezdek's counterexample
\dd\cite{Bezdek} to one of Hsiang's main claims.

Such a practice follows a basic tenet of science about the
dissemination of criticism. This is seen also in the principle that
an author should not referee his own paper: he can never be deemed
objective about his own work. Then similarly he can hardly always
be objective about its correctness. If solely an author should
write, or judge over errata to his paper, is such a rule there to
manifest truth, or authority? In fact, Prof.\ A himself advised
me against complacency as an author, writing in connection to
my experience in \S\ref{S3}: "You should be modest. You believe
your paper is excellent, [$\ldots$] but others may not think so."
This is certainly true, but it makes his attitude towards his own
case appear that much more strange.

There is, however, indeed a difference between critique of published
and unpublished work, which becomes obvious from the journal's point
of view. No journal is responsible for an unpublished manuscript, and
in fact journals seek ways to dispose of submissions. Whether unable or
unwilling, editors are seldom concerned about legitimacy of critiques
toward authors, and generally protect their referees (as seen in the
editorial quote of \S\ref{SI}). This has led certain referees to the
conviction that peer review is a good arena for hitting on "enemies" as
one pleases (instead of, but doing as if discussing mathematical flaws).
The situation changes, though, if a paper is published. On the one
hand, some authors see in a publication by a reputed journal, however
achieved, a certification of supremacy. On the other hand, problems
arising with its publications directly affect the journal. Thus a new
culprit must be sought, who is easily found in the audience. It is not
surprising, therefore, that when authors rebuff an objecting reader,
some editors happily join the chorus. The result is to enforce deciding
how "excellent" one's paper is by whether and where it is published,
while discouraging and obstructing its actual study.

Another issue raised in the M editor's letter to me was length. I
find it a rather questionable argument that form should have priority
over content, and that for length reasons M cannot publish a single
page of corrections, when it publishes hundreds of pages of research
per year. Who is served by preferring long papers to correct ones?
Or are all M's papers ascertainably error-free?

A journal declaring a published paper to have undergone some refereeing
is very far from a guaranty. Even with the purest (and only very
occasionally present) intention to examine a paper mathematically~--
and not politically~-- a referee is not infallible. Often he sees
well his shortcomings and recommends publication relying on authors
and journal to take ultimate responsibility. This makes sense, since
neither can he take such responsibility, when his work is not publicly
disclosed, nor should he, when he acts voluntarily and is not supposed
to receive any credit for a publication. Thus pronouncing refereeing
a stamp of correctness for everybody to believe in is not more than
some editors' absolutory tale.

\section{Common editorial practices\label{cep}}

Of course I have to respect editors of a mathematical class much
higher than my own. Also, we should be inclined to believe that
many of them do their best to make a journal for their readers as
good as they can. It is plausible that many submissions, diverse
referee opinions, a lack of definite evaluation criteria, etc.,
pose difficulties to the publishing selection process of a journal.
On the other hand, an author's life is hardly made easier. A
journal demands exclusive right to handle a manuscript, which is
not to be submitted simultaneously anywhere else, but it guarantees
no time frame for this. Delays to consider a submission (or even
not to consider it, as in my example) are not uncommon. Neither
are comments of a referee written without any degree of politeness.
Even with a positive referee report, or after the author is
requested to submit a revision, the final decision can well be
negative, and a debate over it is deemed useless. At the same
time, the author may need to seek employment, or to compete with
several other researchers working on the same subject.

However allegedly essential, freedom and severity of editorial
action appear to me particularly objectionable in the case
of corrections. One common practice is not to give any serious
mathematical explanation when declining publication. Two types
of argument are often, implicitly or explicitly, appealed to in
such a situation: the limits of printing space, and the journal's
control over quality. However, over corrections such phrases not
only, as discussed in \cite[\S II]{Lang}, prevent meaningful
scientific debate, but they appear questionable even from the
point of view of mere common sense.

If being merciless over submissions is meant to achieve quality,
it is not clear what notion of quality is attained by leaving
errors unfixed (and improperly examined). I hereby clearly separate
quality from bibliometrical rankings or the like \cite{BD,IMU,Ewing},
which seem increasingly perceived as a measurement for it~-- and
which have already documentably led to worrisome editorial conduct
\cite{Mushtaq,Arnold}. Neither is it convincing that errata, which
usually occupy a minor part of the printing space, are mainly
responsible for the pressure of publishing backlog. When journals
(agreeably) declare upon rejecting papers that they cannot publish
everything worthwhile, why do they still prioritize new material
to corrections? Are impact factors at the basis of gradually turning
the reliability of the literature into a liability?

When any scientific reasoning fails, the ultimate argument to
silence all dispute is, as explained, that the editor's (or referee's)
\em{opinion} decides what material is suitable for the journal.
Conveniently, one is then commonly pointed to other journals as a
possibly more appropriate place for one's work (see for this below).
However agreeable or not, such rhetoric might make sense as far as
it concerns new research. But when a journal's opinion is valued
enough to justify it turning away and requesting its own errors to
be addressed somewhere else, I see nothing at all that could not be
justified. What will it start looking like when, according to its
opinion, every journal can feel free to do and let publication of
whatever it pleases? Is the order to the occurring chaos to be
sought throughout the publishing landscape, the internet, or in
someone's private conversations?

The difficulties with publishing in "some independent journal"
are evident from the reactions in \S\ref{S4}. One other journal
quoted a referee calling my erratum directly "a strange paper".
Its purpose was indicated, and its story had been explained to
the editors at submission. Whatever their form or importance,
I do not see what should be so "strange" about corrective efforts.
If a `normal' way to judge scientific material is that the statement
matters, not the proof, is this how one is concerned about what
appeals to the journal's readership? 

I admit that there are recommended, and possibly more appealing
(but far lower than M in political stature) places where I have
personal reasons not to submit my correction. This is not unusual.
The publishing process has now taken essentially every right from
an author, except one: his choice of journal. Even this sole freedom
is constantly targeted by attempts to manipulate authors where to
submit. I say `manipulate', because such efforts are by no means
always in the author's interest. Many of these recommendations
may not be intendedly adversary, but neither are they compellingly
serious and sincere. My experience has shown it \em{vital} for my
publishing activity that I almost \em{continuously disregarded}
such advice.

Of course, this may result in an author's failure to have his
work properly published. In such a case, an attitude taken in the
publishing system, as seen from A's quote in \S\ref{S3}, is that an
author is responsible for securing attention to whatever he does.
Even for a correction, the reply from M directs the question toward
whether ``\em{[I] want} to publish this'' [emphasis added], and away
from whether or what \em{they should} do. Had A, B or M taken proper
action, or would the diverse people aware of the case still convince
them to do so, then submitting my note, and the whole discussion about
who and where to read it, will be unnecessary. However, while these
same editors and referees bounce around a correction like much other
material, few of them seriously entertain "the fundamental problem
of scientists not answering scientific criticisms of their work, not
allowing publication of criticisms, or requiring other scientists
to submit to various authorities" \cite[\S 3]{Lang}\footnote{A
colleague quoted an editor writing to him: "Your paper is unsuitable
for publication because it corrects another paper."}. No matter
how plausible their attitude, the result is to make journals
virtually "clogged" as "ordinary scientific channels [$\ldots$] for
the presentation of scientific challenges". (\em{ibid.})

As indicated, many problems arise also when treating the internet at
large as such a venue. Lacking an appropriate level of reliability,
and deployed as a retreat from publishing responsibility, it will
only deteriorate standards of communicating research. When a
correction is available \em{somewhere} on the internet, can it be
readily located? Will it be permanently maintained? Who stands for
its accuracy? And most importantly, again: will those who published
the mistake actively engage in settling these issues, or will they
believe it~-- and leave it~-- to be someone else's job? Or when "a
couple of experts" are aware of every problem, will they make effort
to honestly and consistently spread the information, and will they
reach anyone who potentially needs it? Not only my case has shown
many "experts" concerned about other things than mathematics.

\section{Conclusion: Privatizing mathematical correctness\label{S8}}

Over a long period Prof.\ A has repeatedly suggested to me that
problems are better quietly entrusted to those in charge, and
that challenges disturb the common atmosphere. There might have
been the needed sense of responsibility around in scientific
community in the climate he grew up 40 years ago to make this a
workable attitude. But I see such responsibility neglected now.
Applying Prof.\ A's principle to his own case, I waited almost 3
years before objecting publicly. What I have seen during that time
is inaction, scientific misrepresentation, and arbitrary power. 

I do not agree that authors who deem me an "enemy" or that I "cannot
afford more enemies" should not seriously discuss their flawed
mathematics. I clarify again that I cannot judge the importance
of my correction, and that it would have been unnecessary, had
mathematical arguments shown otherwise, or had the authors themselves
handled the matter properly. But I do maintain that correctness
issues need open (and not only the authors') evaluation and
jurisdiction. And I do maintain that the proper way to answer
criticism (mathematical and ethical) is not that, "like the video
games[,] one can't shoot fast enough" \cite[p.\ 797]{Lang2} at
(job-seeking) "enemies". I feel not satisfied with scientific
achievement, or the internet, as an exemption from publishing
responsibility. And I do not see covering up flaws in papers as
a good way of polishing journals or research resumes.

Without a job for several years, I have repeatedly found people with
decades of professional experience, including such who advised me
regarding my juvenile behavior, showing looser manners than my
own. There is reason to condemn my indignant reactions in various
situations, as long as some minimum on academic values is respected.
But when I fail to see anything about mathematics to stay out of
all the political intrigue-making, those responsible for maintaining
standards have turned it into an uphill battle. And this problem
neither occurs only due to, nor does it affect only myself.
G.~Perelman, despite being widely respected for his work, retired
from mathematical research, and was quoted regarding his decision
\cite{NasarGruber}: "It is not people who break ethical standards
who are regarded as aliens. It is people like me who are isolated.
[$\ldots$] Of course, there are many mathematicians who are more or
less honest. But almost all of them are conformists. They are more or
less honest, but they tolerate those who are not honest." I consider
this an excellent summary of the state which policies of indulgence,
pleasing and fighting each other in mathematics have led to.

Thus I wrote my account down, in order to lend importance at least
to this question: \emph{Should one tolerate that authors and
journals can evade involvement in correctness issues of their own
publications?} I fear that if this becomes widely adopted, journals
will be there for propaganda, while the introspection of mathematics
will be relegated to the gossip. And mathematicians will be
tempted to politically smuggle papers past scientific control. In
relation to the Baltimore case, P.~Doty \cite{Doty} wrote that
such "attitude towards the responsibility of authors [$\ldots$]
is a critical departure from common standards\dots\ [T]o leave
to others the responsibility of establishing the validity of
what you have published is not only a fundamental retreat from
responsibility but, if it became accepted practice, would erode the
way science works. For [$\ldots$] science moves forward by building
rapidly on what is published on the tentative assumption that it is
correct, not by waiting for others to test each paper's validity."

Indeed, while I cannot seek public attention to every case I consider
annoying, I do not see much decency coming up ahead with (what is
supposed to be called) a "publishing system" looking like this.
Mathematics can only be done by human beings, and may carry their
imperfection. But is it still felt as a common duty, or is it now
a private property? "In this way we risk sliding down toward the
standards of some other professions where {\bf the validity of action
is decided by whether one can get away with it} [boldface added]. For
science to drift toward such a course would be fatal~-- not only to
itself and the inspiration which carries it forward, but to the public
trust which is its provider." (Doty \em{ibid.}\footnote{Lang, who has
extensively quoted Doty's piece, comments on it on \cite[p.339]{Lang2}
that ``we do not `risk' sliding down toward such standards; we have
reached them.''})

\section{Postscriptum: what to do?}

As a response to my letter to the AMS\ \cite{nams}, I had some
discussion with Michael Fried, who pointed me to his article
\cite{Fried}. In order to improve refereeing standards, Fried
proposes a pool of special referees for high-quality journals.
He believes that "behind considerable corruption is a community
neglect of developing significant mathematical skills", which
should be promoted (or preserved) by educating such reviewers.
These referees should be public and receive a small payment. I
insist that "public" should mean disclosing the referee's identity
for a particular paper, and in case the paper is accepted, the
referee can share a minor part of credit and responsibility for
publication. I am strongly convinced that, as long as referees
are kept anonymous, a considerable improvement in their ethics
is not in sight. Thus paying anonymous referees will only aggravate
the problem. I have also some reservations toward a blind-refereeing
approach (where the author is unknown to the referee). This will
pressure authors to write papers with an aim to conceal their own
identity. Such attempts can still often be easily uncovered (in
particular, through the internet), and will only go at the cost
of presentation.

Another point is the need to establish a general administrative
process to deal with issues of correctness (both ethical and
scientific) in mathematics. According to a familiar source, the
AMS has no such channel, and I know of nothing similar at other
mathematical institutions. But my case is not isolated (see, e.g., 
the mentioned articles by T.~Hill), and there exist structures in
other sciences to address such problems, which suggests that these
may deserve more serious attention in mathematics, too. I point out
again, following \cite[\S V.3]{Lang} (as quoted in \S\ref{Sli}), that
civil and academic responsibilities are two different things, which
is why challenging scientific conduct by legal procedures not only
bears enormous risks, but also questionable prospects. The lack of
scientific procedures has created a "legalistic morass" (Lang
\em{ibid.}), where the idea blossoms that as long as "it is hard to
prove it" for the "enemies" what one does, all should be considered
legitimate.

Ultimately, whatever system is devised, I see the central point to
be, as S.\ Lang writes in \cite[Conclusion]{Lang}, that scientists
"uphold the traditional standards of science." And in doing so,
"[t]hey must rely on individual responsibility, and they must create
an atmosphere and conditions under which scientists, both young and
established, can exercise this responsibility without fear~-- fear of
retaliation, fear for their careers, fear for their funding, fear for
their publications, fear of the tensions which come from a challenge,
fear of being uncollegial, whatever. Will they?" (\em{ibid.})


\begin{thebibliography}{20}
% \bibitem[AMS]{M} \emph{American Mathematical Society},
%   \web|http://www.ams.org|

\bibitem[AMS Math\dd{}Rev]{MSN} \emph{American
  Mathematical Society Mathematical Reviews (MathSciNet)},
  \web|http://www.ams.org/publications/math-reviews/math-reviews|

\bibitem[Arnold 2009]{Arnold} Arnold, Douglas N. (2009),
  \emph{Integrity Under Attack: The State of Scholarly Publishing},
  SIAM News, December 4, 2009,
  \webb.|http://www.siam.org/news/news.php?id=1663|

\bibitem[arXiv]{arXiv} \emph{arXiv.org e-Print archive},
  \web|http://arxiv.org/archive/math|

\bibitem[Bernstein]{Bernstein} Bernstein, Dennis S., \emph{A Student's
  Guide to Peer Review},
 \web|http://www.et.byu.edu/~beard/Helps_for_students/peer_review.pdf|

\bibitem[Bernstein2]{Bernstein2} \bysame, \emph{On Review Practice},
  \web|http://aerospace.engin.umich.edu/people/faculty/bernstein/guide/ReviewPractice.pdf|

\bibitem[Bezdek 1997]{Bezdek}
  Bezdek, K. (1997), \emph{Isoperimetric inequalities and the
  dodecahedral conjecture}, Internat. J. Math., Volume 8, Number
  6, 759-780.

\bibitem[Bornmann 2011]{B} Bornmann, L. (2011), \emph{Scientific
  peer review}, Annual Review of Information Science and Technology,
  Volume 45, Chapter 5, \webb.|http://www.lutz-bornmann.de/|

\bibitem[Bornmann-Daniel 2008]{BD} Bornmann, L. and Daniel, H.-D.
  (2008), \emph{What do citation counts measure? A review of studies on
  citing behavior}, Journal of Documentation, Volume 64, Number 1,
  45-80.

\bibitem[Buzzelli 1993]{Bu} Buzzelli, Donald E. (1993), \emph{The
   Definition of Misconduct in Science: A View from NSF}, Science,
   Volume 259, Issue 5095, 584-585, 647-648.

\bibitem[IMU 2008]{IMU} \emph{Citation Statistics: An IMU
  Report} (2008), Notices of the American Mathematical Society
  Volume 55, Number 8, 968-969.

\bibitem[Doty 1991]{Doty} Doty, P. (1991), \emph{Responsibility
  and Weaver et al.}, Nature, Volume 352, 18 July 1991, 183-184.

\bibitem[Ewing 2006]{Ewing} Ewing, John (2006), \emph{Measuring
  Journals}, Notices of the American Mathematical Society, Volume 53,
  Number 9, 1049-1053.

\bibitem[Fried 2007]{Fried} Fried, Michael D. (2007), \emph{Should
  Journals compensate Referees?}, Notices of the AMS, Volume 54,
  Number 6, 585,
  \webb-|http://www.math.uci.edu/~mfried/proplist-ams.html|.

\bibitem[Hales 2006]{Hales} Hales, Thomas C. (2006),
  \emph{Historical Overview of the Kepler Conjecture}, Discrete
  Comput. Geom., Volume 36, 5-20, DOI: 10.1007/s00454-005-1210-2

\bibitem[Hill 2009]{Hill} Hill, Theodore P. (2009),
  \emph{How to Publish Counterexamples in 1 2 3 Easy Steps},
  \webb r|http://www.scribd.com/doc/19819297/How-to-Publish-Counterexamples-in-1-2-3-Easy-Steps|

\bibitem[Hill 2010]{Hill2} \bysame\ (2010), \emph{Hoisting
  the Black Flag}, Letters to the Editor, Notices of the American
  Mathematical Society, Volume 57, Number 1, 7,
  \web|http://www.ams.org/notices/201001/index.html|

\bibitem[Jackson 2002]{Jackson} Jackson, Allyn (2002), \emph{From Preprints
  to E-prints: The Rise of Electronic Preprint Servers in Mathematics},
  Notices of the American Mathematical Society,
  Volume 49, Number 1, 23-32.

\bibitem[Lang 1993]{Lang} Lang, Serge (1993), \emph{Questions of
  Scientific Responsibility: the Baltimore case}, Ethics and Behavior,
  Volume 3, Number 1, 3-72,
  \web|http://www.gatewaycoalition.org/files/Gateway_Project_Moshe_Kam/Resource/DBC.html|

\bibitem[Lang 1998]{Lang2} \bysame\ (1998), \emph{Challenges},
  Springer Verlag, 816 pages.

\bibitem[Mushtaq 2007]{Mushtaq} Mushtaq, Qaiser (2007), \emph{The
  Misuse of the Impact Factor}, Opinion, Notices of the American
  Mathematical Society, Volume 54, Number 7, 821.

\bibitem[Nasar-Gruber 2006]{NasarGruber} Sylvia Nasar and David Gruber
  (2006), \emph{Manifold Destiny. A legendary problem and the battle
  over who solved it}, The New Yorker, August 28, 2006,
  \web|http://www.newyorker.com/archive/2006/08/28/060828fa_fact2|

\bibitem[Odling-Smee et al. 2007]{O+} Odling-Smee, Lucy; Giles, Jim;
  Fuyuno, Ichiko; Cyranoski, David and Marris, Emma (2007), \emph{Where
  are they now?}, Nature, Volume 445, 244-245. DOI: 10.1038/445244a
   
\bibitem[Stoimenow 2010]{nams} Stoimenow, Alexander (2010), {\em
  Honesty in Mathematical Writing}, Letters to the Editor, Notices
  of the AMS, Volume 57, Number 6 (June/July 2010), 703.

\bibitem[Szpiro 2003]{Szpiro} Szpiro, George G. (2003), \emph{Kepler's
  Conjecture: How Some of the Greatest Minds in History Helped
  Solve One of the Oldest Math Problems in the World},
  Hoboken, NJ, John Wiley \& Sons.

\end{thebibliography}
\end{document}